\documentclass[review]{elsarticle}

\usepackage{lineno,hyperref}
\usepackage{amsfonts,amssymb}
\usepackage{mathrsfs}  

\journal{Journal of \LaTeX\ Templates}

\newdefinition{rmk}{Definition}
\newproof{pf}{Proof}
\newproof{pot}{Proof of Theorem \ref{thm2}}

\begin{document}

\begin{frontmatter}

\title{A puzzle on the existence of Lyapunov functions for Limit cycle system
}

\author[firstaddress,secondaryaddress,thirdaddress]{Xiao-Liang Gan}
\ead{ganxiaoliang21@163.com}
\author[firstaddress,secondaryaddress,thirdaddress]{Hao-Yu Wang}
\author[secondaryaddress,thirdaddress,fourthaddress]{Ping Ao\corref{mycorrespondingauthor}}
\cortext[mycorrespondingauthor]{Corresponding author}
\ead{aoping@sjtu.edu.cn}

\address[firstaddress]{Department of Mathematics, Shanghai University, Shanghai 200444, China.}

\address[secondaryaddress]{Institute of Systems Science, Shanghai University, Shanghai 200444, China.}

\address[thirdaddress]{Shanghai Center for Quantitative Life Sciences, Shanghai University, Shanghai 200444, China.}

\address[fourthaddress]{Department of Physics, Shanghai University, Shanghai 200444, China.}

\begin{abstract}
Although the limit cycle have been studied for more than 100 years, the existence of its Lyapunov function is still poorly understood.
By considering a common limit cycle system, a puzzle related to the existence of Lyapunov functions for Limit cycle system is proposed and studied in this paper: The divergence is not equal to zero, while the trajectory can be infinite loop on limit cycle.
It will be discussed from three aspects.
Firstly, the definition of Lyapunov function is concerned, and a new version of definition is given which is equivalent to the usual one and better adapts to the properties of Lyapunov function.
Secondly, two criteria(divergence and dissipation power) are discussed from the perspective of dissipation, and it reaches that they are not consistent in representing the dissipation.
Thirdly, by studying the motion of charged massless particle in electromagnetic field, it obtains that the limit cycle is an isopotential line on which the charged particles can move infinitely.
Such discussions above may provide an understanding on the puzzle about the existence of Lyapunov functions for limit cycle system.

\end{abstract}


\end{frontmatter}


\section{Introduce a puzzle}
Consider a common limit cycle system\cite{Perko-2001-p195,ZhangandFeng-2000-73,Chicone-2006-p96}
\begin{equation}\label{1k1}
\left\{ {\begin{array}{*{20}{l}}
{{{\dot x}_1} = {f_1}({x_1},{x_2}) =  - {x_2} + {x_1}[1 - (x_1^2 + x_2^2)]}\\
{{{\dot x}_2} = {f_2}({x_1},{x_2}) = {x_1} + {x_2}[1 - (x_1^2 + x_2^2)]}
\end{array}} \right..
\end{equation}
By the polar transformation$\left\{ {\begin{array}{*{20}{l}}
{{x_1} = r\cos \theta }\\
{{x_2} = r\sin \theta }
\end{array}} \right.,$ the (\ref{1k1}) can be translated into
\begin{equation}\label{1k6}
\left\{ \begin{array}{l}
\frac{{dr}}{{dt}} = r - {r^3}\\
\frac{{d\theta }}{{dt}} =   1
\end{array} \right..
\end{equation}
Obviously, it has a limit cycle ${x_1^2 + x_2^2}= 1$.

The divergence on the limit cycle of system(\ref{1k1}) is
\begin{eqnarray}\label{1k2}
{\left. {\nabla  \cdot f} \right|_{{x_1^2 + x_2^2}= 1}} &=&{\left. {\left[ {\frac{{\partial {f_1}}}{{\partial {x_1}}} + \frac{{\partial {f_2}}}{{\partial {x_2}}}} \right]} \right|_{x_1^2 + x_2^2 = 1}} \nonumber\\
   &=& {\left. {\left[ {(1 - 3x_1^2 - x_2^2) + (1 - x_1^2 - 3x_2^2)} \right]} \right|_{x_1^2 + x_2^2 = 1}} \nonumber\\
   &=&  {\left. {\left[ {2(1 - 2x_1^2 - 2x_2^2)} \right]} \right|_{x_1^2 + x_2^2 = 1}}  =  - 2 \ne 0.
\end{eqnarray}

A puzzle arises: The divergence is not equal to zero, while the trajectory can be infinite loop on limit cycle.
For this phenomenon, researchers have two problems:
\begin{description}
  \item[Question 1:] Whether the system(\ref{1k1}) has a Lyapunov function, or not?
  \item[Question 2:] If the Lyapunov function exists, how to use it to analysis and understand the infinite motion on limit cycle?
\end{description}
Some researchers have explored these two questions and gaven some opinions.
For the Question 1, such as, Strogatz\cite{Strogatz-2014-p203} believes that for a system, if the Lyapunov function exists, there is no closed orbit.
Hirsch, Smale and Devaney\cite{HirschSmalDevaney-2012-p225} confirm that if the planar system has a strict Lyapunov function for a planar system, then there are no limit cycles.
Teschl\cite{Teschl-2011-p109} has the similar view with Hirsch.
In a word, the Lyapunov function they mentioned is just for fixed points\cite{Wiggins-2003-p22,Hsu-2013-p140,AlligoodSauerYorke-1996-p305}.
For the Question 2, Zhu and Cai\cite{ZhuCai-2017-p284} think that although the limit cycle is also a periodic orbital, the periodic solution and the limit cycle here have different meanings.
Hirsch, Smale and Devaney\cite{HirschSmalDevaney-2012-p225} point a property of limit cycle that when $t\rightarrow +\infty $ the limit cycle is a limit set. Not all closed orbits have this property, such as the center of linear system.
Lefschetz\cite{Lefschetz-1962-p223} points out that a limit-cycle is a closed path which is not a member of a continuous family of closed paths, i.e. an isolated closed path.
However, they gave no specific explanation to the reason why orbitals move infinitely on limit cycle.
To find an appropriate answers to questions 1 and 2 is difficult, is there an appropriate strategy or approach to address the above problems?  The answer is positive.
Ao\cite{Ao-2004-p25,YuanAo-2012-p7010,YuanMaYuanAo-2014-p10505} divided the dynamic system $\dot x = f(x)$ into three parts from the perspective of mechanics:
\begin{equation}\label{1k3}
(S + T)\dot x =  - \nabla \phi,
\end{equation}
here, $x = {\left[ {{x_1}, \cdots ,{x_n}} \right]^\tau }$, $\tau$ is the transpose symbol, the friction matrix $S$ corresponds to dissipation, the transverse matrix $T$ corresponds to Lorentz force and the potential function $\phi$ is equivalent to Lyapunov function.
Based on the new construction(\ref{1k3}), Yuan et al.\cite{YuanMaYuanAo-2014-p10505} have answered the first question and find a Lyapunov function for system(\ref{1k1})
\begin{equation}\label{1k4}
\phi ({x_1},{x_2}) = \frac{1}{4}(x_1^2 + x_2^2)(x_1^2 + x_2^2 - 2).
\end{equation}
The Figure 1 and Figure 2 are the trajectory picture in vector field and Lyapunov function of (\ref{1k1}), respectively.
\begin{figure}[htbp]
\centering\begin{minipage}[t]{55mm}
\includegraphics[width=55mm]{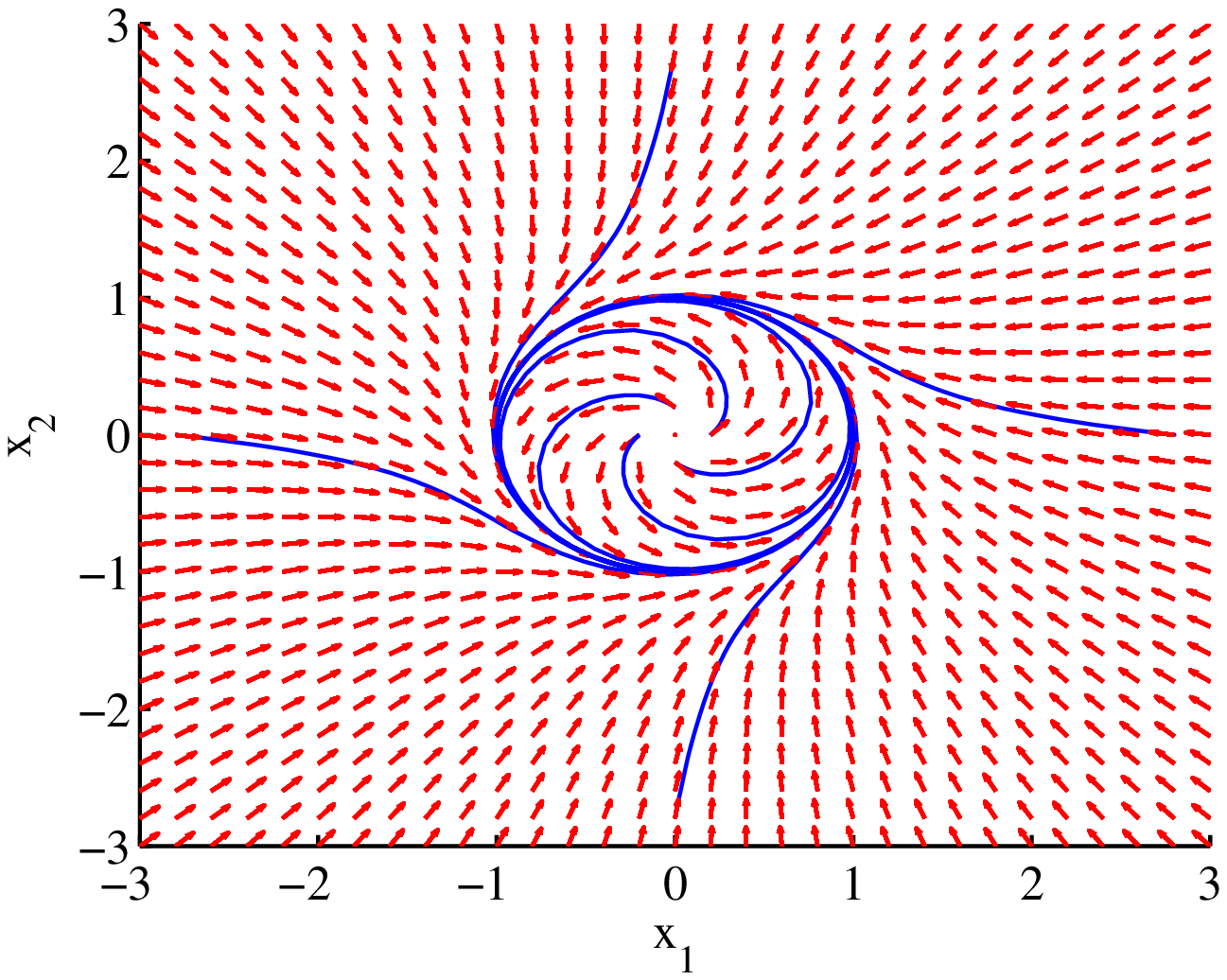}
\caption{Vector fields and trajectory of(\ref{1k1})}
\label{figp1}
\end{minipage}
\qquad
\begin{minipage}[t]{55mm}
\includegraphics[width=55mm]{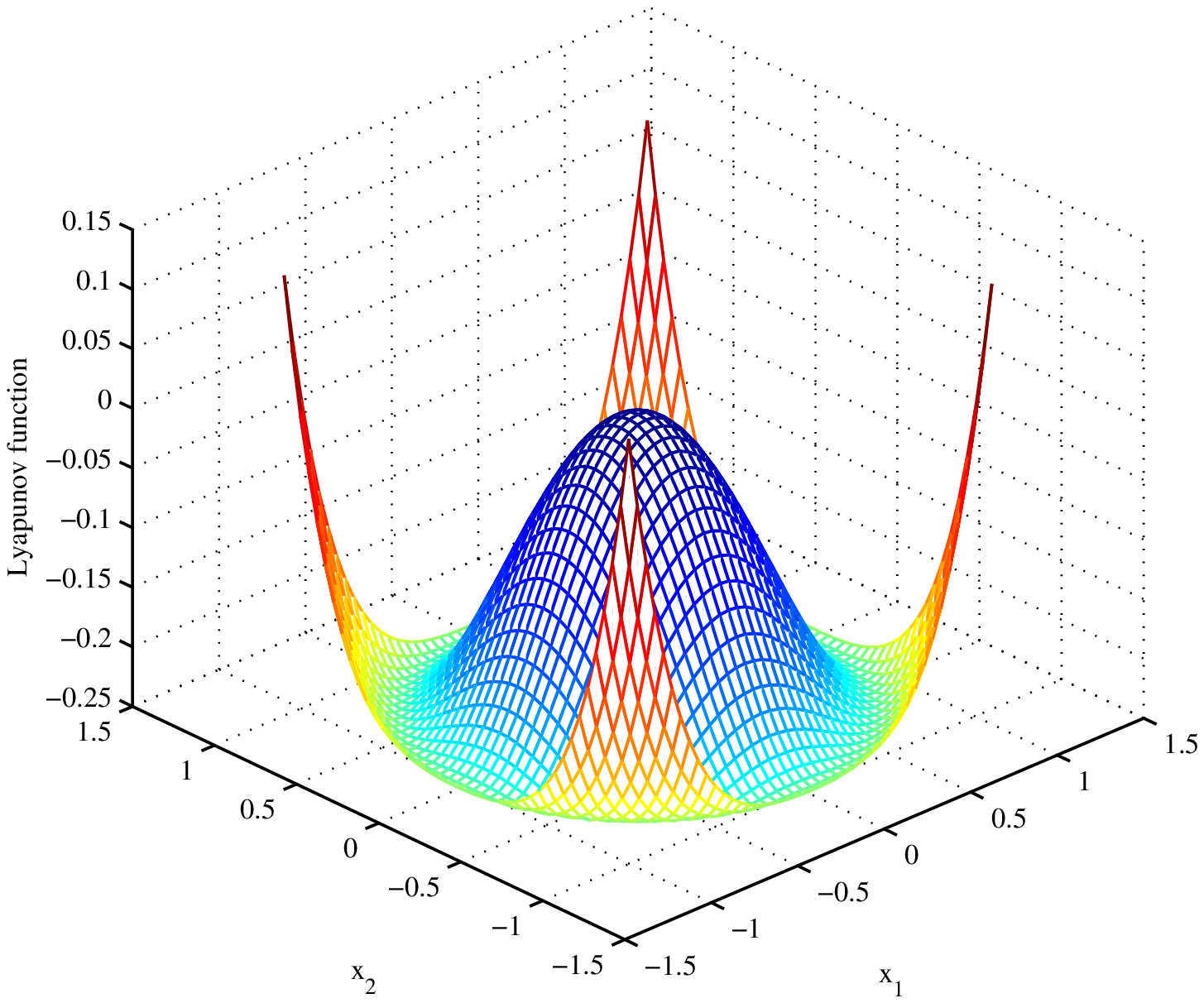}
\caption{Lyapunov function of(\ref{1k1})}
\label{figp2}
\end{minipage}
\end{figure}

Furthermore, we can verify that the Lyapunov function is decreasing by following trajectory
\begin{equation}\label{1k5}
\frac{{d\phi ({x_1},{x_2})}}{{dt}} =  - (x_1^2 + x_2^2){(x_1^2 + x_2^2 - 1)^2} \le 0.
\end{equation}

In the following, combining with the new construction(\ref{1k3}) and results of Yuan et al.\cite{YuanMaYuanAo-2014-p10505}, we will continue to analyze and explain the puzzle from three aspects in section 2. The conclusion is in section 3.
\section{Analyze and explain the puzzle}
\subsection{From the definition of the Lyapunov function}
Firstly, the definition of Lyapunov function in most textbooks is given.

\textbf{Definition 1\cite{Wiggins-2003-p22,Hsu-2013-p140,AlligoodSauerYorke-1996-p305}:} Consider the following vector field
\begin{equation}\label{1k7}
\dot x = f(x),x \in \mathbb{R}^{n}.
\end{equation}
Let ${x^ * }$ be a fixed point of (\ref{1k7}) and let $V:\Omega \to \mathbb{R}$ be a ${C^1}$ function defined on some neighborhood $\Omega$ of ${x^ * }$ such that
\begin{description}
  \item[(i)] $V({x^ * }) = 0$ and $V(x) > 0$ if $x \ne {x^ * }$.
  \item[(ii)] $\dot V(x) \le 0$ in $\Omega - \{ {x^ * }\} $.
\end{description}
Then $x^ *$ is stable, and the $V(x)$ is called a Lyapunov function for $x^ *$.

However, there are three kinds of dynamic systems: fixed point, limit cycle and chaos.
Wolfram\cite{Wolfram-2002-p961} classified the dynamical systems in this way, too.
The Lyapunov function $V(x)$ just can be used for the qualitative analysis of fixed point. So the promotion of Lyapunov function has a deeper and broader practical significance.

Here, only the fixed point and limit cycle are considered, and chaos is not discussed.

\textbf{Definition 2}
For a smooth autonomous system $\dot x = f(x)$, let ${x^ {**} }$ be a limit set which is a fixed point or a limit cycle. And $\phi(x)$ is a continuous differentiable function from the state space to $\mathbb{R}$.  If it satisfies
\begin{description}
  \item[(1)] for all $x $ in the state space, the Lyapunov function does not increase along the trajectories, that is
  \begin{equation}\label{1j1017}
  \dot \phi (x) = \frac{{d\phi }}{{dt}} \le 0;
  \end{equation}
  \item[(2)] for all ${x}$, $\phi (x)$ has an infimum.
\end{description}
Then, $\phi(x)$ is the generalized Lyapunov function.


\begin{flushleft}
\textbf{Remark}\quad  From the above generalized definition, it can be seen that the fixed point and limit cycle are generalized as the limit set.
The condition (1) ensures that it doesn't increase along trajectory;
The condition (2) generalizes the definition of positive definite.
This definition of Lyapunov function is logically consistent with the one in dynamical system textbooks.
\end{flushleft}

\subsection{Analyze the cause of the puzzle from the perspective of system dissipation}
It's easy to understand that divergence does not equal zero for the Lyapunov function of fixed point.
However, when the limit set is limit cycle, it is hard to imagine.
Combined with the decomposition (\ref{1k3}), we will analyze the cause from the perspective of system dissipation.

There are two criteria of dissipation:

\textbf{Criterion 1} Dissipation is defined as dissipation power\cite{GoldsteinPooleSafko-2001-p24}:
\begin{eqnarray}\label{1j1004}
{H_P}&=& {F_{friction}} \cdot ( - \dot x) \nonumber\\
     &=& ( - S\dot x) \cdot ( - \dot x)  \nonumber\\
     &=& {\dot x^\tau }S\dot x,
\end{eqnarray}
here, $S$ is a positive semidefinite friction matrix. 
This method is often used in physics.

\textbf{Criterion 2} Dissipation is defined as divergence\cite{Huang-2010-p35}:
\begin{eqnarray}\label{1j1005}
div f(x)&=& \sum {\frac{{\partial {f_i}}}{{\partial {x_i}}}},
\end{eqnarray}
this method is often used in mathematics.

\textbf{Remark}\quad
Since the equivalence between potential function and Lyapunov function has been proved by Yuan et al.\cite{YuanMaYuanAo-2014-p10505}, we can further verify that the rate of change of the Lyapunov function$\left| {\frac{{d\phi }}{{dt}}} \right|$ is equal to dissipation power$H_p$. That is
\begin{equation}\label{1k8}
\left| {\frac{{d\phi }}{{dt}}} \right|=H_p.
\end{equation}

In the following, just the planar dynamical system is taken into consideration, and the divergence and dissipated power will be derived out and compared to see whether they are consistent in representing the dissipation or not.


\begin{flushleft}
\textbf{2.2.1} $\dot x = f(x)$ is a linear system
\end{flushleft}

The two dimension linear dynamical system can be expressed as
\begin{equation}\label{1j526}
\dot x = f(x) = Ax, A~ is~ a~ constant~ matrix, x = {\left[ {{x_1},{x_2}} \right]^\tau }.
\end{equation}

Kwon, Ao and Thouless\cite{KwonAoThouless-2005-p13029} have discussed the construction of Lyapunov function of linear system. And some necessary formulas are given
\begin{equation}\label{1j526}
A=  - \left[ {D + Q} \right]U =  - {\left[ {S + T} \right]^{ - 1}}U,
\end{equation}
\begin{equation}\label{1j721}
{[D + Q]^{ - 1}} = S + T,
\end{equation}
\begin{equation}\label{1j528}
AQ + Q{A^\tau } = AD - D{A^\tau }.
\end{equation}
here, $Q$ is an antisymmetric matrix, $D$ is a semidefinite diffusion matrix, $S$ semidefinite friction matrix,
$T$ is an antisymmetric transverse matrix, $U$ is a symmetric matrix. And only the $Q$ is unknown, the $S$ and $T$ can be derived by $D$ and $Q$.

By (\ref{1j526})and (\ref{1j721}), it derives
\begin{equation}\label{1j527}
U =  - {\left[ {D + Q} \right]^{ - 1}}A =  - {A^\tau }{\left[ {D - Q} \right]^{ - 1}}.
\end{equation}

Set $A = \left[ {\begin{array}{*{20}{c}}
{{a_{11}}}&{{a_{12}}}\\
{{a_{21}}}&{{a_{22}}}
\end{array}} \right],Q = \left[ {\begin{array}{*{20}{c}}
0&{{Q_{12}}}\\
{ - {Q_{12}}}&0
\end{array}} \right],D = \left[ {\begin{array}{*{20}{c}}
{{d_{11}}}&{{d_{12}}}\\
{{d_{12}}}&{{d_{22}}}
\end{array}} \right]$,
and it satisfies
${d_{11}},{d_{22}} \ge 0,{d_{11}}{d_{22}} - d_{12}^2 \ge 0.$

By(\ref{1j528}), it has
\begin{equation}\label{1j603}
({a_{11}} + {a_{22}}) {Q_{12}}=  - {a_{21}}{d_{11}} + ({a_{11}} - {a_{22}}){d_{12}} + {a_{12}}{d_{22}},
\end{equation}
here, only $Q_{12}$ is unknown, ${a_{11}},{a_{12}},{a_{21}},{a_{22}}$ are known, and ${d_{11}},{d_{12}},{d_{22}} $ satisfy
${d_{11}},{d_{22}} \ge 0,{d_{11}}{d_{22}} - d_{12}^2 \ge 0$.
Then, it has
\begin{equation}\label{1j616}
\left\{ \begin{array}{l}
({a_{11}} + {a_{22}}){Q_{12}} =  - {a_{21}}{d_{11}} + ({a_{11}} - {a_{22}}){d_{12}} + {a_{12}}{d_{22}}\\
{d_{11}},{d_{22}} \ge 0,{d_{11}}{d_{22}} - d_{12}^2 \ge 0
\end{array} \right..
\end{equation}

The matrix $A = \left[ {\begin{array}{*{20}{c}}
{{a_{11}}}&{{a_{12}}}\\
{{a_{21}}}&{{a_{22}}}
\end{array}} \right]$  has the following four types of Jordan's normal form\cite{MaZhou-2013-p100}
\begin{equation}\label{1j615}
\left[ {\begin{array}{*{20}{c}}
{{\lambda _{\rm{1}}}}&{\rm{0}}\\
{\rm{0}}&{{\lambda _{\rm{2}}}}
\end{array}} \right],\left[ {\begin{array}{*{20}{c}}
{{\lambda _{\rm{1}}}}&{\rm{0}}\\
{\rm{0}}&{{\lambda _{\rm{1}}}}
\end{array}} \right],\left[ {\begin{array}{*{20}{c}}
{{\lambda _{\rm{1}}}}&{\rm{0}}\\
1&{{\lambda _{\rm{1}}}}
\end{array}} \right],\left[ {\begin{array}{*{20}{c}}
\alpha &\beta \\
{ - \beta }&\alpha
\end{array}} \right],
\end{equation}
here, $\lambda_1, \lambda_2$,$\alpha  \pm \beta i$ are the different eigenvalues of $A$. $\alpha  \pm \beta i$ is a pair of conjugate eigenvalues, and $\beta \ne 0$.

In the following, let matrix $A $ be one of the forms (\ref{1j615}), and get the corresponding divergence $div f(x)$ and dissipation power ${H_P} $.
\begin{description}
  \item[(1)] When $A = \left[ {\begin{array}{*{20}{c}}
   {{\lambda _{\rm{1}}}}&{\rm{0}}\\
   {\rm{0}}&{{\lambda _{\rm{2}}}}
   \end{array}} \right]$, (\ref{1j616}) can be rewritten as
  \begin{equation}\label{1j617}
  \left\{ \begin{array}{l}
  ({\lambda _{\rm{1}}} + {\lambda _2}){Q_{12}} = ({\lambda _{\rm{1}}} - {\lambda _2}){d_{12}}\\
  {d_{11}},{d_{22}} \ge 0,{d_{11}}{d_{22}} - d_{12}^2 \ge 0
  \end{array} \right.,
  \end{equation}
  \begin{description}
    \item[(i)]  If ${\lambda _1} + {\lambda _2} \ne 0$, it can get $Q_{12}$ by (\ref{1j617})
    \begin{equation}\label{1j618}
    {Q_{12}} = \frac{{{\lambda _{\rm{1}}} - {\lambda _2}}}{{{\lambda _{\rm{1}}} + {\lambda _2}}}{d_{12}},
    \end{equation}
    then, it has
    \begin{eqnarray}\label{1h1}
      Q &=& \frac{{{\lambda _{\rm{1}}} - {\lambda _2}}}{{{\lambda _{\rm{1}}} + {\lambda _2}}}{d_{12}}\left[ {\begin{array}{*{20}{c}}
0&1\\
{ - 1}&0
\end{array}} \right],  \\
    D + Q &=& \left[ {\begin{array}{*{20}{c}}
{{d_{11}}}&{\frac{{2{\lambda _{\rm{1}}}}}{{{\lambda _{\rm{1}}} + {\lambda _2}}}{d_{12}}}\\
{\frac{{2{\lambda _2}}}{{{\lambda _{\rm{1}}} + {\lambda _2}}}{d_{12}}}&{{d_{22}}}
\end{array}} \right], \\
     {\left[ {D + Q} \right]^{ - 1}}  &=& \frac{{{{({\lambda _{\rm{1}}} + {\lambda _2})}^2}}}{{{d_{11}}{d_{22}}{{({\lambda _{\rm{1}}} + {\lambda _2})}^2} - 4{\lambda _{\rm{1}}}{\lambda _2}d_{12}^2}}\left[ {\begin{array}{*{20}{c}}
{{d_{22}}}&{ - \frac{{2{\lambda _{\rm{1}}}}}{{{\lambda _{\rm{1}}} + {\lambda _2}}}{d_{12}}}\\
{ - \frac{{2{\lambda _2}}}{{{\lambda _{\rm{1}}} + {\lambda _2}}}{d_{12}}}&{{d_{11}}}
\end{array}} \right] \\
   S &=& \frac{{{{\left[ {D + Q} \right]}^{ - 1}} + {{\left\{ {{{\left[ {D + Q} \right]}^{ - 1}}} \right\}}^\tau }}}{2}  \nonumber \\
     &=&  \frac{{{{({\lambda _{\rm{1}}} + {\lambda _2})}^2}}}{{{d_{11}}{d_{22}}{{({\lambda _{\rm{1}}} + {\lambda _2})}^2} - 4{\lambda _{\rm{1}}}{\lambda _2}d_{12}^2}}\left[ {\begin{array}{*{20}{c}}
    {{d_{22}}}&{ - {d_{12}}}\\
    { - {d_{12}}}&{{d_{11}}}
     \end{array}} \right] ,
    \end{eqnarray}
\begin{eqnarray}\label{1j624}
U &=&  - \frac{{{{({\lambda _{\rm{1}}} + {\lambda _2})}^2}}}{{{{({\lambda _{\rm{1}}} + {\lambda _2})}^2}{d_{11}}{d_{22}} - 4{\lambda _{\rm{1}}}{\lambda _2}d_{12}^2}}\left[ {\begin{array}{*{20}{c}}
{{\lambda _{\rm{1}}}{d_{22}}}&{ - \frac{{2{\lambda _{\rm{1}}}{\lambda _2}}}{{{\lambda _{\rm{1}}} + {\lambda _2}}}{d_{12}}}\\
{ - \frac{{2{\lambda _{\rm{1}}}{\lambda _2}}}{{{\lambda _{\rm{1}}} + {\lambda _2}}}{d_{12}}}&{{\lambda _2}{d_{11}}}
\end{array}} \right].
\end{eqnarray}
  By  $\nabla \phi (x) = Ux$, it can derive the Lyapunov function
    \begin{equation}\label{1j625}
     \phi  =  - \frac{{{{({\lambda _{\rm{1}}} + {\lambda _2})}^2}({d_{22}}{\lambda _{\rm{1}}}x_1^2 - \frac{{4{\lambda _{\rm{1}}}{\lambda _2}}}{{{\lambda _{\rm{1}}} + {\lambda _2}}}{d_{12}}{x_1}{x_2} + {d_{11}}{\lambda _2}x_2^2)}}{{{d_{11}}{d_{22}}{{({\lambda _{\rm{1}}} + {\lambda _2})}^2} - 4{\lambda _{\rm{1}}}{\lambda _2}d_{12}^2}},
    \end{equation}
    Then, it can verify
    \begin{eqnarray}\label{1j626}
     \frac{{d\phi }}{{dt}} &=&\frac{{\partial \phi }}{{\partial {x_1}}}{{\dot x}_1} + \frac{{\partial \phi }}{{\partial {x_2}}}{{\dot x}_2} \nonumber\\
       &=& - \frac{{{{({\lambda _{\rm{1}}} + {\lambda _2})}^2}({d_{22}}\lambda _1^2x_1^2 - 2{\lambda _{\rm{1}}}{\lambda _2}{d_{12}}{x_1}{x_2} + {d_{11}}\lambda _2^2x_2^2)}}{{{{({\lambda _{\rm{1}}} + {\lambda _2})}^2}{d_{11}}{d_{22}} - 4{\lambda _1}{\lambda _2}d_{12}^2}}\le 0.
    \end{eqnarray}
   Then, the corresponding divergence $div f(x)$ and dissipation power${H_P} $ are derived
  \begin{eqnarray}\label{1h9}
  {H_P} &=&  {{\dot x}^\tau }S\dot x = \frac{{{{({\lambda _{\rm{1}}} + {\lambda _2})}^2}({d_{22}}\lambda _1^2x_1^2 - 2{d_{12}}{\lambda _1}{\lambda _2}{x_1}{x_2} + {d_{11}}\lambda _2^2x_2^2)}}{{{d_{11}}{d_{22}}{{({\lambda _{\rm{1}}} + {\lambda _2})}^2} - 4{\lambda _{\rm{1}}}{\lambda _2}d_{12}^2}}\ne 0,\\
    div f(x)  &=& \nabla  \cdot f = trace (A)={{\lambda _{\rm{1}}} + {\lambda _2}} \ne 0.
  \end{eqnarray}
  And then, it can verify
   \begin{equation}\label{1k9}
   \left| {\frac{{d\phi }}{{dt}}} \right|=H_p.
   \end{equation}
     In this case, $div f(x)$ and ${H_P}$ are not zero at the same time. And they all show that the system is dissipative. So these two criteria are consistent in representing the dissipation of the system.
    \item[(ii)] If ${\lambda _{\rm{1}}} + {\lambda _2}=0$, it can get $Q_{12}$ by (\ref{1j617})
    \begin{equation}\label{1j724}
    d_{12}=0, Q_{12} ~is~ an~ arbitrary~ real~ number.
    \end{equation}
  then, it has
  \begin{eqnarray}\label{1h2}
    D &=& \left[ {\begin{array}{*{20}{c}}
{{d_{11}}}&0\\
0&{{d_{22}}}
\end{array}} \right],{d_{11}},{d_{22}}\ge 0 , \\
   Q  &=& \left[ {\begin{array}{*{20}{c}}
0&{{Q_{12}}}\\
{ - {Q_{12}}}&0
\end{array}} \right], \\
   D + Q &=& \left[ {\begin{array}{*{20}{c}}
{{d_{11}}}&{{Q_{12}}}\\
{ - {Q_{12}}}&{{d_{22}}}
\end{array}} \right], \\
   {\left[ {D + Q} \right]^{ - 1}} &=& \frac{1}{{{d_{11}}{d_{22}} + Q_{12}^2}} \left[ {\begin{array}{*{20}{c}}
{{d_{22}}}&{ - {Q_{12}}}\\
{{Q_{12}}}&{{d_{11}}}
\end{array}} \right], \\
   S &=&  \frac{1}{{{d_{11}}{d_{22}} + Q_{12}^2}}\left[ {\begin{array}{*{20}{c}}
{{d_{22}}}&0\\
0&{{d_{11}}}
\end{array}} \right],
  \end{eqnarray}
  \begin{equation}\label{1j631}
    U =  - {\left[ {D + Q} \right]^{ - 1}}A =  - \frac{{{\lambda _1}}}{{{d_{11}}{d_{22}} + Q_{12}^2}}\left[ {\begin{array}{*{20}{c}}
{{d_{22}}}&{{Q_{12}}}\\
{{Q_{12}}}&{ - {d_{11}}}
\end{array}} \right],
    \end{equation}
   By  $\nabla \phi (x) = Ux$, it can derive the Lyapunov function
    \begin{equation}\label{1j632}
    \phi  =  - \frac{{{\lambda _1}({d_{22}}x_1^2 + 2{Q_{12}}{x_1}{x_2} - {d_{11}}x_2^2)}}{{{d_{11}}{d_{22}} + Q_{12}^2}},
    \end{equation}
   Then, it can verify
    \begin{equation}\label{1J633}
    \frac{{d\phi }}{{dt}} = \frac{{\partial \phi }}{{\partial {x_1}}}{{\dot x}_1} + \frac{{\partial \phi }}{{\partial {x_2}}}{{\dot x}_2} =  - \frac{{\lambda _1^2({d_{22}}x_1^2 + {d_{11}}x_2^2)}}{{{d_{11}}{d_{22}} + Q_{12}^2}} \le 0.
    \end{equation}
  The corresponding divergence $div f(x)$ and dissipation power${H_P} $ are derived
  \begin{eqnarray}\label{1h10}
  {H_P} &=& {{\dot x}^\tau }S\dot x =\frac{{\lambda _1^2({d_{22}}x_1^2 + {d_{11}}x_2^2)}}{{{d_{11}}{d_{22}} + Q_{12}^2}} \ne 0,  \\
   div f(x) &=& \nabla  \cdot f = trace (A)={{\lambda _{\rm{1}}} + {\lambda _2}}=0.
  \end{eqnarray}
  \end{description}
  Then, it can verify
   \begin{equation}\label{1k10}
   \left| {\frac{{d\phi }}{{dt}}} \right|=H_p.
   \end{equation}
  In this case, $div f(q)=0$ doesn't mean ${H_P}=0$. So these two criteria are not consistent in representing the dissipation of the system.
  \item[(2)] When$A = \left[ {\begin{array}{*{20}{c}}
   {{\lambda _{\rm{1}}}}&{\rm{0}}\\
   {\rm{0}}&{{\lambda _{\rm{1}}}}
   \end{array}} \right]$, (\ref{1j616}) can be rewritten as
   \begin{equation}\label{1j640}
   \left\{ {\begin{array}{*{20}{l}}
{2{\lambda _1}{Q_{12}} = 0}\\
{{d_{11}},{d_{22}} \ge 0,{d_{11}}{d_{22}} - d_{12}^2 \ge 0}
\end{array}} \right.
   \end{equation}
   \begin{description}
     \item[(i)] If$\lambda_1 \ne 0$, it can get $Q_{12}$ by (\ref{1j640})
   \begin{equation}\label{1j641}
     {Q_{12}} = 0,
   \end{equation}
  then, it has
  \begin{eqnarray}\label{1h3}
   Q  &=& 0,  \\
  D + Q &=&  \left[ {\begin{array}{*{20}{c}}
{{d_{11}}}&{{d_{12}}}\\
{{d_{12}}}&{{d_{22}}}
\end{array}} \right], \\
  {\left[ {D + Q} \right]^{ - 1}}&=& \frac{1}{{{d_{11}}{d_{22}} - d_{12}^2}}\left[ {\begin{array}{*{20}{c}}
{{d_{22}}}&{ - {d_{12}}}\\
{ - {d_{12}}}&{{d_{11}}}
\end{array}} \right], \\
    S &=&\frac{1}{{{d_{11}}{d_{22}} - d_{12}^2}}\left[ {\begin{array}{*{20}{c}}
{{d_{22}}}&{ - {d_{12}}}\\
{ - {d_{12}}}&{{d_{11}}}
\end{array}} \right],
  \end{eqnarray}
   \begin{equation}\label{1j647}
   U =  - {\left[ {D + Q} \right]^{ - 1}}A =  - \frac{{{\lambda _1}}}{{{d_{11}}{d_{22}} + d_{12}^2}}\left[ {\begin{array}{*{20}{c}}
{{d_{22}}}&{ - {d_{12}}}\\
{ - {d_{12}}}&{{d_{11}}}
\end{array}} \right],
   \end{equation}
   By $\nabla \phi (x) = Ux$, it can derive the Lyapunov function
   \begin{equation}\label{1j648}
   \phi  =  - \frac{{{\lambda _1}}}{{{d_{11}}{d_{22}} + d_{12}^2}}({d_{22}}x_1^2 - 2{d_{12}}{x_1}{x_2} + {d_{11}}x_2^2),
   \end{equation}
  Then, it can verify
   \begin{equation}\label{1j649}
   \frac{{d\phi }}{{dt}} = \frac{{\partial \phi }}{{\partial {x_1}}}{{\dot x}_1} + \frac{{\partial \phi }}{{\partial {x_2}}}{{\dot x}_2} =  - \frac{{\lambda _1^2({d_{22}}x_1^2 - 2{d_{12}}{x_1}{x_2} + {d_{11}}x_2^2)}}{{{d_{11}}{d_{22}} + d_{12}^2}} \le 0.
   \end{equation}
  The corresponding divergence $div f(x)$ and dissipation power${H_P} $ are derived
  \begin{eqnarray}\label{1h11}
   {H_P} &=& {{\dot x}^\tau }S\dot x = \frac{{\lambda _1^2({d_{22}}x_1^2 - 2{d_{12}}{x_1}{x_2} + {d_{11}}x_2^2)}}{{{d_{11}}{d_{22}} - d_{12}^2}}\ne 0 . \\
   div f(x)&=&  \nabla  \cdot f = trace (A)=2{\lambda _{\rm{1}}} \ne 0 .
  \end{eqnarray}
  Then, it can verify
   \begin{equation}\label{1k11}
   \left| {\frac{{d\phi }}{{dt}}} \right|=H_p.
   \end{equation}
  In this case, $div f(x)$ and ${H_P}$ are not zero at the same time. And they all show that the system is dissipative. So these two criteria are consistent in representing the dissipation of the system.
     \item[(ii)] If $\lambda_1 =0$, the system is conservative. It's easy to know
  \begin{eqnarray}\label{1h12}
  {H_P} &=& {{\dot x}^\tau }S\dot x \equiv 0,  \\
    div f(x)  &=& \nabla  \cdot f = trace (A)=2{\lambda _{\rm{1}}} \equiv 0 .
  \end{eqnarray}
  Then, it can verify
   \begin{equation}\label{1k12}
   \left| {\frac{{d\phi }}{{dt}}} \right|=H_p.
   \end{equation}
    In this case, they are equal to zero at the same time. $div f(x)\equiv 0$ and ${H_P}\equiv0$ are consistent in representing the system being conservative.
   \end{description}
  \item[(3)]When $A = \left[ {\begin{array}{*{20}{c}}
   {{\lambda _{\rm{1}}}}&{\rm{0}}\\
   {\rm{1}}&{{\lambda _{\rm{1}}}}
   \end{array}} \right]$, (\ref{1j616}) can be rewritten as
   \begin{equation}\label{1j652}
   \left\{ {\begin{array}{*{20}{l}}
{2{\lambda _1}{Q_{12}} =  - {d_{11}}}\\
{{d_{11}},{d_{22}} \ge 0,{d_{11}}{d_{22}} - d_{12}^2 \ge 0}
\end{array}} \right..
   \end{equation}
   \begin{description}
     \item[(i)] If $\lambda_1 \ne 0$,  it can get $Q_{12}$ by (\ref{1j652})
     \begin{equation}\label{1j653}
     {Q_{12}} = \frac{{ - {d_{11}}}}{{2{\lambda _1}}},
     \end{equation}
  then, it has
  \begin{eqnarray}\label{1h4}
   Q  &=& \frac{{ - {d_{11}}}}{{2{\lambda _1}}}\left[ {\begin{array}{*{20}{c}}
0&1\\
{ - 1}&0
\end{array}} \right], \\
  D + Q &=& \left[ {\begin{array}{*{20}{c}}
{{d_{11}}}&{{d_{12}} - \frac{{{d_{11}}}}{{2{\lambda _1}}}}\\
{{d_{12}} + \frac{{{d_{11}}}}{{2{\lambda _1}}}}&{{d_{22}}}
\end{array}} \right], \\
    {\left[ {D + Q} \right]^{ - 1}}   &=& \frac{{4\lambda _1^2}}{{4\lambda _1^2{d_{11}}{d_{22}} - 4\lambda _1^2d_{12}^2 + d_{11}^2}}\left[ {\begin{array}{*{20}{c}}
{{d_{22}}}&{ - {d_{12}} + \frac{{{d_{11}}}}{{2{\lambda _1}}}}\\
{ - {d_{12}} - \frac{{{d_{11}}}}{{2{\lambda _1}}}}&{{d_{11}}}
\end{array}} \right], \\
   S  &=& \frac{{4\lambda _1^2}}{{4\lambda _1^2{d_{11}}{d_{22}} - 4\lambda _1^2d_{12}^2 + d_{11}^2}}\left[ {\begin{array}{*{20}{c}}
{{d_{22}}}&{ - {d_{12}}}\\
{ - {d_{12}}}&{{d_{11}}}
\end{array}} \right],
  \end{eqnarray}
   \begin{eqnarray}\label{1659}
    U &=& - \frac{{4\lambda _1^2}}{{4\lambda _1^2{d_{11}}{d_{22}} - 4\lambda _1^2d_{12}^2 + d_{11}^2}}\left[ {\begin{array}{*{20}{c}}
{{d_{22}}{\lambda _{\rm{1}}} - {d_{12}} + \frac{{{d_{11}}}}{{2{\lambda _1}}}}&{ - {d_{12}}{\lambda _1} + \frac{{{d_{11}}}}{2}}\\
{ - {d_{12}}{\lambda _1} + \frac{{{d_{11}}}}{2}}&{{d_{11}}{\lambda _1}}
\end{array}} \right],
   \end{eqnarray}
   By $\nabla \phi (x) = Ux$, it can derive the Lyapunov function
     \begin{equation}\label{1j660}
    \phi  =  - \frac{{4\lambda _1^2\left[ {({d_{22}}{\lambda _{\rm{1}}} - {d_{12}} + \frac{{{d_{11}}}}{{2{\lambda _1}}})x_1^2 + 2( - {d_{12}}{\lambda _1} + \frac{{{d_{11}}}}{2}){x_1}{x_2} + {d_{11}}{\lambda _1}x_2^2} \right]}}{{4\lambda _1^2{d_{11}}{d_{22}} - 4\lambda _1^2d_{12}^2 + d_{11}^2}},
     \end{equation}
  Then, it gets
     \begin{equation}\label{1j661}
     \frac{{d\phi }}{{dt}} =- \frac{{4\lambda _1^2\left[ {({d_{22}}\lambda _1^2 - 2{d_{12}}{\lambda _{\rm{1}}} + {d_{11}})x_1^2 + 2( - {d_{12}}{\lambda _1} + {d_{11}}){\lambda _1}{x_1}{x_2} + {d_{11}}\lambda _1^2x_2^2} \right]}}{{4\lambda _1^2{d_{11}}{d_{22}} - 4\lambda _1^2d_{12}^2 + d_{11}^2}},
     \end{equation}
      the root discriminant of ${({d_{22}}\lambda _1^2 - 2{d_{12}}{\lambda _{\rm{1}}} + {d_{11}})x_1^2 + 2( - {d_{12}}{\lambda _1} + {d_{11}}){\lambda _1}{x_1}{x_2} + {d_{11}}\lambda _1^2x_2^2}=0$ is
     \begin{eqnarray}\label{1j662}
      \Delta &=& 4{( - {d_{12}}{\lambda _1} + {d_{11}})^2}\lambda _1^2 - 4({d_{22}}\lambda _1^2 - 2{d_{12}}{\lambda _{\rm{1}}} + {d_{11}}){d_{11}}\lambda _1^2 \nonumber \\
        &=&  4\lambda _1^4\left[ {d_{12}^2 - {d_{11}}{d_{22}}} \right] \le 0,
     \end{eqnarray}
     it has ${({d_{22}}\lambda _1^2 - 2{d_{12}}{\lambda _{\rm{1}}} + {d_{11}})x_1^2 + 2( - {d_{12}}{\lambda _1} + {d_{11}}){\lambda _1}{x_1}{x_2} + {d_{11}}\lambda _1^2x_2^2} \geq 0$, then it obtains
     \begin{equation}\label{1j663}
     \frac{{d\phi }}{{dt}}=\frac{{\partial \phi }}{{\partial {x_1}}}{{\dot x}_1} + \frac{{\partial \phi }}{{\partial {x_2}}}{{\dot x}_2}\leq 0.
     \end{equation}
  The corresponding divergence $div f(q)$ and dissipation power${H_P} $ are derived
  \begin{eqnarray}\label{1h15}
  {H_P}  &=&  \frac{{4\lambda _1^2[(\lambda _1^2d_{22}^2 - 2{\lambda _1}{d_{12}} + {d_{11}})x_1^2 + 2( - \lambda _1^2{d_{12}} + {\lambda _1}{d_{11}}){x_1}{x_2} +
      \lambda _1^2{d_{11}}x_2^2]}}{{4\lambda _1^2{d_{11}}{d_{22}} - 4\lambda _1^2d_{12}^2 + d_{11}^2}} \nonumber\\
      &\neq& 0, \\
   div f(x)  &=& \nabla  \cdot f = trace (A)=2{\lambda _{\rm{1}}} \ne 0 .
  \end{eqnarray}
   Then, it can verify
   \begin{equation}\label{1k13}
   \left| {\frac{{d\phi }}{{dt}}} \right|=H_p.
   \end{equation}
  In this case, $div f(x)$ and ${H_P}$ are not zero at the same time. And they all show that the system is dissipative. So these two criteria are consistent in representing the dissipation of the system.
     \item[(ii)] If $\lambda_1 = 0$, it can get $Q_{12}$ by (\ref{1j652})
     \begin{equation}\label{1j666}
     d_{11}=d_{12}=0, Q_{12}~is~ an~ arbitrary~ real~ number.
     \end{equation}
  then, it has
  \begin{eqnarray}\label{1h5}
    D &=& \left[ {\begin{array}{*{20}{c}}
0&0\\
0&{{d_{22}}}
\end{array}} \right],  \\
   Q &=& \left[ {\begin{array}{*{20}{c}}
0&{{Q_{12}}}\\
{ - {Q_{12}}}&0
\end{array}} \right],   \\
   D + Q &=& \left[ {\begin{array}{*{20}{c}}
0&{{Q_{12}}}\\
{ - {Q_{12}}}&{{d_{22}}}
\end{array}} \right],  \\
  {\left[ {D + Q} \right]^{ - 1}}&=&  \frac{1}{{Q_{12}^2}}\left[ {\begin{array}{*{20}{c}}
{{d_{22}}}&{ - {Q_{12}}}\\
{{Q_{12}}}&0
\end{array}} \right], \\
   S &=& \frac{1}{{Q_{12}^2}}\left[ {\begin{array}{*{20}{c}}
{{d_{22}}}&0\\
0&0
\end{array}} \right],
  \end{eqnarray}
  \begin{equation}\label{1j673}
     U =  - {\left[ {D + Q} \right]^{ - 1}}A = \frac{1}{{{Q_{12}}}}\left[ {\begin{array}{*{20}{c}}
1&0\\
0&0
\end{array}} \right],
     \end{equation}
   By $\nabla \phi (x) = Ux$, it can derive the Lyapunov function
     \begin{equation}\label{1j674}
     \phi  = \frac{1}{{{Q_{12}}}}x_1^2,
     \end{equation}
  Then, it can verify that the Lyapunov function is decreasing by trajectory
     \begin{equation}\label{1j675}
     \frac{{d\phi }}{{dt}} = \frac{{\partial \phi }}{{\partial {x_1}}}{{\dot x}_1} + \frac{{\partial \phi }}{{\partial {x_2}}}{{\dot x}_2} = \frac{{2{x_1}{{\dot x}_1}}}{{{Q_{12}}}}= \frac{{2{\lambda _1}x_1^2}}{{{Q_{12}}}} = 0,
     \end{equation}
  then, the corresponding divergence $div f(x)$ and dissipation power${H_P} $ are derived
  \begin{eqnarray}\label{1h13}
  {H_P}  &=&  {{\dot x}^\tau }S\dot x=  {x^\tau }{A^\tau }SAx \nonumber \\
   &=& \frac{1}{{Q_{12}^2}}\left[ {\begin{array}{*{20}{c}}
{{x_1}}&{{x_2}}
\end{array}} \right]\left[ {\begin{array}{*{20}{c}}
0&1\\
0&0
\end{array}} \right]\left[ {\begin{array}{*{20}{c}}
{{d_{22}}}&0\\
0&0
\end{array}} \right]\left[ {\begin{array}{*{20}{c}}
0&0\\
1&0
\end{array}} \right]\left[ {\begin{array}{*{20}{c}}
{{q_1}}\\
{{q_2}}
\end{array}} \right] \nonumber \\
&\equiv& 0, \\
   div f(x)  &=& \nabla  \cdot f = trace (A)=2{\lambda _{\rm{1}}} \equiv 0 .
  \end{eqnarray}
      Then, it can verify
   \begin{equation}\label{1k14}
   \left| {\frac{{d\phi }}{{dt}}} \right|=H_p.
   \end{equation}
  In this case, they are equal to zero at the same time. $div f(x)\equiv 0$ and ${H_P}\equiv0$ are consistent in representing the system being conservative.
   \end{description}
  \item[(4)]When $A=\left[ {\begin{array}{*{20}{c}}
\alpha &\beta \\
{ - \beta }&\alpha
\end{array}} \right]$, (\ref{1j616}) can be rewritten as
\begin{equation}\label{1j678}
\left\{ {\begin{array}{*{20}{l}}
{2\alpha {Q_{12}} = \beta ({d_{11}} + {d_{22}})}\\
{{d_{11}},{d_{22}} \ge 0,{d_{11}}{d_{22}} - d_{12}^2 \ge 0,\beta  \ne 0}
\end{array}} \right..
\end{equation}
\begin{description}
  \item[(i)] If $\alpha \ne 0$, it can get $Q_{12}$ by (\ref{1j678})
  \begin{equation}\label{1j679}
  {Q_{12}} = \frac{{\beta ({d_{11}} + {d_{22}})}}{{2\alpha }},
  \end{equation}
  then, it has
  \begin{eqnarray}\label{1h6}
   Q &=&  \frac{{\beta ({d_{11}} + {d_{22}})}}{{2\alpha }}\left[ {\begin{array}{*{20}{c}}
0&1\\
{ - 1}&0
\end{array}} \right], \\
   D + Q &=& \left[ {\begin{array}{*{20}{c}}
{{d_{11}}}&{{d_{12}} + \frac{{\beta ({d_{11}} + {d_{22}})}}{{2\alpha }}}\\
{{d_{12}} - \frac{{\beta ({d_{11}} + {d_{22}})}}{{2\alpha }}}&{{d_{22}}}
\end{array}} \right],  \\
   {\left[ {D + Q} \right]^{ - 1}}  &=&  \frac{{4{\alpha ^2}}}{{4{\alpha ^2}({d_{11}}{d_{22}} - d_{12}^2) + {\beta ^2}{{({d_{11}} + {d_{22}})}^2}}}\left[ {\begin{array}{*{20}{c}}
{{d_{22}}}&{ - {d_{12}} - \frac{{\beta ({d_{11}} + {d_{22}})}}{{2\alpha }}}\\
{ - {d_{12}} + \frac{{\beta ({d_{11}} + {d_{22}})}}{{2\alpha }}}&{{d_{11}}}
\end{array}} \right], \\
   S &=& \frac{{4{\alpha ^2}}}{{4{\alpha ^2}({d_{11}}{d_{22}} - d_{12}^2) + {\beta ^2}{{({d_{11}} + {d_{22}})}^2}}}\left[ {\begin{array}{*{20}{c}}
{{d_{22}}}&{ - {d_{12}}}\\
{ - {d_{12}}}&{{d_{11}}}
\end{array}} \right],
  \end{eqnarray}
  \begin{scriptsize}
  \begin{eqnarray}\label{1j685}
   U &=& \frac{{ - 4{\alpha ^2}}}{{4{\alpha ^2}({d_{11}}{d_{22}} - d_{12}^2) + {\beta ^2}{{({d_{11}} + {d_{22}})}^2}}}\left[ {\begin{array}{*{20}{c}}
{{d_{22}}\alpha  + {d_{12}}\beta  + \frac{{{\beta ^2}({d_{11}} + {d_{22}})}}{{2\alpha }}}&{ - {d_{12}}\alpha  + \frac{{\beta ({d_{22}} - {d_{11}})}}{2}}\\
{ - {d_{12}}\alpha  + \frac{{\beta ({d_{22}} - {d_{11}})}}{2}}&{{d_{11}}\alpha  - {d_{12}}\beta  + \frac{{{\beta ^2}({d_{11}} + {d_{22}})}}{{2\alpha }}}
\end{array}} \right],
  \end{eqnarray}
  \end{scriptsize}
   By $\nabla \phi (x) = Ux$, it can derive the Lyapunov function
     \begin{scriptsize}
     \begin{equation}\label{1j686}
     \phi  = \frac{{ - 4{\alpha ^2}\left\{ {[{d_{22}}\alpha  + {d_{12}}\beta  + \frac{{{\beta ^2}({d_{11}} + {d_{22}})}}{{2\alpha }}]x_1^2 + 2[ - {d_{12}}\alpha  + \frac{{\beta ({d_{22}} - {d_{11}})}}{2}]{x_1}{x_2} + [{d_{11}}\alpha  - {d_{12}}\beta  + \frac{{{\beta ^2}({d_{11}} + {d_{22}})}}{{2\alpha }}]x_2^2} \right\}}}{{4{\alpha ^2}({d_{11}}{d_{22}} - d_{12}^2) + {\beta ^2}{{({d_{11}} + {d_{22}})}^2}}},
     \end{equation}
     \end{scriptsize}
     then, it gets
     \begin{tiny}
     \begin{eqnarray}\label{1j687}
     \frac{{d\phi }}{{dt}}&=&\frac{{ - 4{\alpha ^2}\left\{ {({d_{22}}{\alpha ^2} + 2{d_{12}}\alpha \beta  + {d_{11}}{\beta ^2})x_1^2 + 2[{d_{12}}({\beta ^2} - {\alpha ^2}) + \alpha \beta ({d_{22}} - {d_{11}})]{x_1}{x_2} + ({d_{11}}{\alpha ^2} + {d_{22}}{\beta ^2} - 2{d_{12}}\alpha \beta )x_2^2} \right\}}}{{4{\alpha ^2}({d_{11}}{d_{22}} - d_{12}^2) + {\beta ^2}{{({d_{11}} + {d_{22}})}^2}}},
     \end{eqnarray}
     \end{tiny}
     the root discriminant of equation $\left\{  *  \right\} = 0$ corresponding to the formula in (\ref{1j687}) is
     \begin{tiny}
     \begin{eqnarray}\label{1j688}
      \Delta  &=& 4{[{d_{12}}({\beta ^2} - {\alpha ^2}) + \alpha \beta ({d_{22}} - {d_{11}})]^2} - ({d_{22}}{\alpha ^2} + 2{d_{12}}\alpha \beta  + {d_{11}}{\beta ^2})({d_{11}}{\alpha ^2} + {d_{22}}{\beta ^2} - 2{d_{12}}\alpha \beta )\nonumber  \\
        &=& 4{({\alpha ^2} + {\beta ^2})^2}(d_{12}^2 - {d_{11}}{d_{22}}) \le 0,
     \end{eqnarray}
     \end{tiny}
     then, it has
     \begin{equation}\label{1j689}
     \frac{{d\phi }}{{dt}} = \frac{{\partial \phi }}{{\partial {x_1}}}{{\dot x}_1} + \frac{{\partial \phi }}{{\partial {x_2}}}{{\dot x}_2} \le 0,
     \end{equation}
  then, the corresponding divergence $div f(x)$ and dissipation power${H_P} $ are derived
     \begin{tiny}
    \begin{eqnarray}\label{1j723}
     {H_P} &=& {{\dot x}^\tau }S\dot x  \nonumber\\
       &=& \frac{{4{\alpha ^2}\left\{ {({\alpha ^2}{d_{22}} + 2\alpha \beta {d_{12}} + {\beta ^2}{d_{11}})x_1^2 + 2[\alpha \beta {d_{22}} - ({\alpha ^2} - {\beta ^2}){d_{12}} - \alpha \beta {d_{11}}]{x_1}{x_2} + ({\beta ^2}{d_{22}} - 2\alpha \beta {d_{12}} + {\alpha ^2}{d_{11}})x_2^2} \right\}}}{{4{\alpha ^2}({d_{11}}{d_{22}} - d_{12}^2) + {\beta ^2}{{({d_{11}} + {d_{22}})}^2}}} \nonumber\\
       &\ne& 0  , \\
     div f(x) &=& \nabla  \cdot f = trace (A)=2\alpha \ne 0  .
    \end{eqnarray}
    \end{tiny}
    Then, it can verify
   \begin{equation}\label{1k15}
   \left| {\frac{{d\phi }}{{dt}}} \right|=H_p.
   \end{equation}
   In this case, $div f(x)$ and ${H_P}$ are not zero at the same time. And they all show that the system is dissipative. So these two criteria are consistent in representing the dissipation of the system.
  \item[(ii)] If $\alpha=0$, it can get $Q_{12}$ by (\ref{1j678})
  \begin{equation}\label{1j692}
  {Q_{12}}~is~ an~ arbitrary~ real~ number,{d_{11}} = {d_{22}} = {d_{12}} = 0,
  \end{equation}
  then, it has
  \begin{eqnarray}\label{1h7}
  Q &=& \left[ {\begin{array}{*{20}{c}}
0&{{Q_{12}}}\\
{ - {Q_{12}}}&0
\end{array}} \right],\\
   D  &=& 0, \\
  D + Q &=& \left[ {\begin{array}{*{20}{c}}
0&{{Q_{12}}}\\
{ - {Q_{12}}}&0
\end{array}} \right], \\
   {[D + Q]^{ - 1}} &=& \frac{1}{{Q_{12}^2}}\left[ {\begin{array}{*{20}{c}}
0&{ - {Q_{12}}}\\
{{Q_{12}}}&0
\end{array}} \right], \\
   S  &=& \frac{{{{\left[ {D + Q} \right]}^{ - 1}} + {{\left\{ {{{\left[ {D + Q} \right]}^{ - 1}}} \right\}}^\tau }}}{2} \equiv 0,
  \end{eqnarray}
   \begin{equation}\label{1j698}
  U =  - {\left[ {D + Q} \right]^{ - 1}}A =  - \frac{\beta }{{{Q_{12}}}}\left[ {\begin{array}{*{20}{c}}
1&0\\
0&1
\end{array}} \right],
  \end{equation}
  By $\nabla \phi (x) = Ux$, it can derive the Lyapunov function
  \begin{equation}\label{1j699}
  \phi  =  - \frac{\beta }{{{Q_{12}}}}(x_1^2 + x_2^2),
  \end{equation}
   The Lyapunov function is decreasing by following trajectory
  \begin{eqnarray}\label{1j700}
   \frac{{d\phi }}{{dt}}&=& \frac{{\partial \phi }}{{\partial {x_1}}}{{\dot x}_1} + \frac{{\partial \phi }}{{\partial {x_2}}}{{\dot x}_2} =  - \frac{{2\beta }}{{{Q_{12}}}}({x_1}{{\dot x}_1} + {x_2}{{\dot x}_2})\nonumber \\
     &=& - \frac{{2\beta }}{{{Q_{12}}}}[{x_1}\beta {x_2} + {x_2}( - \beta {x_1})] = 0,
  \end{eqnarray}
  then, the corresponding divergence $div f(x)$ and dissipation power${H_P} $ are derived
  \begin{eqnarray}\label{1h14}
   {H_P} &=&  {{\dot x}^\tau }S\dot x \equiv 0. \\
   div f(x) &=& \nabla  \cdot f = trace (A)\equiv0  .
  \end{eqnarray}
  Then, it can verify
   \begin{equation}\label{1k16}
   \left| {\frac{{d\phi }}{{dt}}} \right|=H_p.
   \end{equation}
  In this case, they are equal to zero at the same time. $div f(x)\equiv 0$ and ${H_P}\equiv0$ are consistent in representing the system being conservative.
\end{description}
\end{description}
In a word, when the fixed point is a saddle point, the divergence$div f(x)=0$ while the dissipation power ${H_P}\ne 0$.
So these two criteria are not always consistent in representing the dissipation of the linear system.


\begin{flushleft}
\textbf{2.2.2} $\dot x = f(x)$ is a non-linear system
\end{flushleft}

In this case, we will analyze it in combination with system (\ref{1k1}).
By the paper \cite{YuanMaYuanAo-2014-p10505}, $S(x)$ can be constructed
if the Lyapunov function $\phi(x) $is known for the nonlinear system
\begin{equation}\label{1j7}
S(x) =  - \frac{{\nabla \phi  \cdot f}}{{f \cdot f}}I,
\end{equation}
here, $I$ represents the identity matrix.

Combined with (\ref{1k4}) and (\ref{1j7}), it can obtain the $S(x)$ of system (\ref{1k1})
\begin{equation}\label{1j722}
S = \frac{{{{(1 - x_1^2 - x_2^2)}^2}}}{{1 + {{(1 - x_1^2 - x_2^2)}^2}}}\left[ {\begin{array}{*{20}{c}}
1&0\\
0&1
\end{array}} \right].
\end{equation}
Furthermore,
\begin{eqnarray}\label{1h16}
{H_P}   &=& (x_1^2 + x_2^2){(x_1^2 + x_2^2 - 1)^2},  \\
div f&=& \nabla  \cdot f = \frac{{\partial {f_1}}}{{\partial {x_1}}} + \frac{{\partial {f_2}}}{{\partial {x_2}}} = 2(1 - 2x_1^2 - 2x_2^2).
\end{eqnarray}
Then, combined with (\ref{1k5})and (\ref{1h16}), it can verify
\begin{equation}\label{1k16}
\left| {\frac{{d\phi }}{{dt}}} \right|=H_p.
\end{equation}
We can find that on limit cycle $x_1^2 + x_2^2 = 1$, it gets ${\left. {\nabla  \cdot f} \right|_{x_1^2 + x_2^2 = 1}} = -2 \ne 0$, $ {H_P} =0$.
So these two criteria are not consistent in representing the dissipation of the nonlinear dynamical systems.

To sum up, these two criteria are not consistent in representing the dissipation.

\subsection{ Analyze the movement on the limit cycle from the physical point of view}
In the following, we will discuss and analyze the movement on the limit cycle from the physical point of view.

Yuan et al.\cite{YuanMaYuanAo-2014-p10505} have given a physical point of view on dynamical system. For the sake of completeness, the process is given, and we will use it to analyze and explain the motion on the limit cycle.

For a charged massless particle in the electromagnetic field, its motion Newton equation is
\begin{equation}\label{1j500}
{F_{driving}} = m\ddot x = 0,
\end{equation}
here, ${F_{driving}}$ is usually divided into two parts: friction force and conservative force
\begin{equation}\label{1j501}
{F_{driving}} = {F_{conservative}} + {F_{friction}}.
\end{equation}

The $F_{friction}$ corresponds to the dissipative part in (\ref{1j501}), and it is expressed as
\begin{equation}\label{1j502}
{F_{friction}} =  - S\dot x,
\end{equation}
 $S$ is a positive semi-definite matrix which ensures the property that friction impedes the motion of an object. So it's called the friction matrix.

The conservative force is derived by lorentz force $ e\dot x \times B $ and electric field force $- \nabla \phi (x)$ which is the potential gradient
\begin{equation}\label{1j503}
{F_{conservative}} = e\dot x \times B + [ - \nabla \phi (x)],
\end{equation}
here, $\nabla  \times x = {\partial _i}{x_j} - {\partial _j}{x_i}$.

By (\ref{1j501}), it can get
\begin{equation}\label{1j505}
S\dot x + eB \times \dot x= - \nabla \phi (x) .
\end{equation}

Change the lorentz force $B \times \dot x$ to $T\dot x$, $T$ is an anti-symmetric matrix and named Lorentz-Ao matrix. Then (\ref{1j505}) can be represented as
\begin{equation}\label{1j506}
 [S(x) + T(x)]\dot x =- \nabla \phi (x) .
\end{equation}
There is another form of (\ref{1j506}) in paper \cite{YuanAo-2012-p7010}
\begin{equation}\label{1j507}
\dot x =  - [D(x) + Q(x)]\nabla \phi (x).
\end{equation}
here, $D(x)$ is a semidefinite diffusion matrix, $Q(x)$ an anti-symmetric matrix which denotes Poisson bracket.


Combine with \cite{YuanMaYuanAo-2014-p10505} and the Lyapunov function (\ref{1k4}), it can derive the corresponding $S, T, D, Q$
\begin{equation}\label{1j511}
S = \frac{{{{(1 - x_1^2 - x_2^2)}^2}}}{{1 + {{(1 - x_1^2 - x_2^2)}^2}}}\left[ {\begin{array}{*{20}{c}}
1&0\\
0&1
\end{array}} \right],
\end{equation}
here, the friction matrix $S$ corresponds to dissipation.
\begin{equation}\label{1j512}
T = \frac{{(1 - x_1^2 - x_2^2)}}{{1 + {{(1 - x_1^2 - x_2^2)}^2}}}\left[ {\begin{array}{*{20}{c}}
0&1\\
{ - 1}&0
\end{array}} \right],
\end{equation}
here, the Lorentz-Ao matrix corresponds to the action of the magnetic field.
\begin{equation}\label{1j513}
 - \nabla \phi  =  - \left[ {\begin{array}{*{20}{c}}
{\frac{{\partial \nabla \phi }}{{\partial x_1}}}\\
{\frac{{\partial \nabla \phi }}{{\partial x_2}}}
\end{array}} \right] = \left[ {\begin{array}{*{20}{c}}
{{x_1}(1 - x_1^2 - x_2^2)}\\
{{x_2}(1 - x_1^2 - x_2^2)}
\end{array}} \right],
\end{equation}
here, the potential gradient $ - \nabla \phi $ corresponds to the action of the electric field.
\begin{equation}\label{1j523}
D(x) =   \left[ {\begin{array}{*{20}{c}}
1&0\\
0&1
\end{array}} \right],
\end{equation}
here, $D(x)$ is the diffusion matrix.
\begin{equation}\label{1j524}
Q(x) = \frac{1}{{1 - x_1^2 - x_2^2}}\left[ {\begin{array}{*{20}{c}}
0&1\\
{ - 1}&0
\end{array}} \right].
\end{equation}
here, $Q(x)$ corresponds to the Poisson bracket.


Here, analyze the motion of charged particles on limit cycles
\begin{itemize}
  \item The friction part:${F_{friction}} =  - S\dot x=0$ , so ${H_P} = {{\dot x}^\tau }S\dot x = 0$ ;
  \item The conservative force part:
  \begin{footnotesize}
  \begin{itemize}
    \item on the limit cycle, the electric potential is equal, so the electric field forces do no work.
    \item Lorentz forces are everywhere perpendicular to the direction of motion of the particle and do no work.
  \end{itemize}
  \end{footnotesize}
  To sum up, the system is conservative on limit cycle, and the charged particle in electromagnetic field can move along the isopotential line.
  Here, the isopotential line is the limit cycle, so it moves in circle.

  \textbf{Remark:} By(\ref{1j506}), the magnetic field intensity is proportional to the electric field intensity on the limit cycle: $T(x)\dot x =  - \nabla \phi (x)$.
   It only affects the size of the circle of motion, not changes the energy of the system.
\end{itemize}

\section{Conclusions}
Based on the decomposition of the dynamical system of Ao, we analyze the puzzle about the existence of Lyapunov functions for limit cycle system.
The generalized definition of Lyapunov function is logically consistent with the definition in dynamical systems textbooks, and it better adapts to the properties of Lyapunov function. By discussing the two criteria of dissipation(divergence and dissipation power), we obtain that they are not consistent in representing the dissipation. The movement on limit cycle is analyzed from the physical point of view, and it obtains that the system is conservative on the limit cycle and the limit cycle is an isopotential line, so it can be infinite motion on limit cycle.
What's more, it may provide an understanding on the puzzle about the existence of Lyapunov function for limit cycle system.

\section*{Acknowledgements}
Thanks members in the Institute of Systems Science of Shanghai University for discussion, particularly, with Xin-Jian Xu.

\section*{References}

\bibliography{mybibfiletwo}

\end{document}